\documentclass[11pt, oneside, a4paper]{article}
\usepackage{amsmath,amssymb,amscd}
\usepackage{latexsym,amsfonts}
\usepackage{xypic}
\usepackage[all]{xy}
\newtheorem{theorem}{Theorem}
\newtheorem{lemma}[theorem]{Lemma}
\newtheorem{cor}[theorem]{Corollary}
\newtheorem{prop}[theorem]{Proposition}
\newtheorem{defn}[theorem]{Definition}
\newtheorem{ex}[theorem]{Example}
\newtheorem{rem}[theorem]{Remark}

\newcommand{\supp}{\mathop{\mathrm{supp}}}
\newcommand{\Fi}{{\mathbb F}}
\newcommand{\Z}{{\mathbb Z}}
\newcommand{\N}{{\mathbb N}}
\DeclareMathOperator{\Aut}{Aut}
\DeclareMathOperator{\Sym}{Sym}
\DeclareMathOperator{\Ker}{ker}
\DeclareMathOperator{\acl}{acl}
\DeclareMathOperator{\dcl}{dcl}
\DeclareMathOperator{\cl}{cl}
\DeclareMathOperator{\Imm}{im}
\DeclareMathOperator{\Th}{Th}
\newcommand{\tp}{\textup{tp}}
\def\ind{\mathop{\mathpalette\Ind{}}}

\newcommand{\cvd}{\hspace*{\fill}
    {\rm \hbox{\vrule height 0.2 cm width 0.2cm}}}

\def\Ind{\setbox0=\hbox{$x$}\kern\wd0\hbox to 0pt{\hss$\mid$\hss}
\lower.9\ht0\hbox to 0pt{\hss$\smile$\hss}\kern\wd0}

\title{Failure of $n$-uniqueness: a family of examples}

\author{Elisabetta Pastori \\ 
Dipartimento di Matematica,\\
Universit\`a degli Studi di Torino, \\
Via Carlo Alberto, 10 10123\\
Torino, Italy \\
{\sf elisabetta.pastori@unito.it}\\
Pablo Spiga \\
Centre for Mathematics of Symmetry and Computation\\
The University of Western Australia\\
35 Stirling Highway\\
Crawley WA 6009\\
Australia\\
{\sf spiga@maths.uwa.edu.au}}

\begin{document}
\date{}

\maketitle

\begin{abstract}
In this paper, the connections between model theory and the
theory of infinite permutation groups (see~\cite{Au}) are used to study the
$n$-existence and the $n$-uniqueness for  
$n$-amalgamation problems of stable
theories. We show that, for any
$n\geq 2$, there exists a  stable theory having $(k+1)$-existence
and $k$-uniqueness, for every $k\leq n$, but has neither
$(n+2)$-existence nor $(n+1)$-uniqueness. 
In particular, this generalizes the example, for $n=2$, due to
E.Hrushovski  given
in~\cite{PKM}. 
\end{abstract}

\maketitle

\section{Introduction}\label{sec0}

Considerable work (e.g. \cite{BK}, \cite{PKM}, \cite{EO}, \cite{Hr},
\cite{She}) has explored higher amalgamation properties for stable and
simple theories.  In this paper we analyze uniqueness and existence
properties for a countable family of stable theories. In contrast to
previous methods our  approach uses
group-theoretic techniques. We begin by giving some basic
definitions. 

Let $T$ be a complete and simple $L$-theory with quantifier
elimination. We denote by $\mathcal{C}_T$ the category of algebraically
closed substructures of models of $T$ with embeddings as
morphisms. Also, given $n\in \mathbb{N}$, we denote by  $P(n)$ the
partially ordered 
set of all subsets of $\{1, \dots ,n\}$ and by
$P(n)^-$ the set $P(n)\setminus \{1,\dots,n\}$.  
 
 An \emph{$n$-amalgamation problem} over
 $\acl(\emptyset)$ is a functor $a:P(n)^-\rightarrow \mathcal{C}_T$ such
 that
\begin{description}
\item[$(i)$] $a(\emptyset)=\acl(\emptyset)$;
\item[$(ii)$] whenever $s_1,s_2,s_3\in P(n)^-$ and $(s_1\cap
  s_2)\subset s_3$, the algebraically closed sets
  $a(s_1),a(s_2)$ are independent over $a(s_1\cap s_2)$ within $a(s_3)$; 
\item[$(iii)$] $a(s)=\acl\{a(i)\, |\, i\in s\}$, for every $s\in P(n)^-$.
\end{description}
In here we denote by $\acl(A)$ the algebraic closure of $A$ in
$T^{\textrm{eq}}$. We recall that the objects of $P(n)^{-}$ (viewed as
a category) are simply the elements of $P(n)^{-}$. Also, the
morphisms of $P(n)^{-}$ are the inclusions $\iota_{s,t}:s\hookrightarrow t$,
for every  $s,t\in P(n)^{-}$ with $s\subseteq t$. In particular, an
$n$-amalgamation problem assigns a morphism $$a_{s,t}:a(s)\to a(t),$$
to every $s,t\in P(n)^{-}$ with  $s\subseteq t$. The morphism
$a_{s,t}$ is called \emph{transition map} and, by functoriality, we have
$$a_{s_2,s_3}\circ a_{s_1,s_2}=a_{s_1,s_3},$$ for every $s_1,s_2,s_3\in
P(n)^{-}$ with $s_1\subseteq s_2\subseteq s_3$. By definition, the morphisms
in $\mathcal{C}_T$ are the embeddings, that is, 
$a_{s,t}$ is the restriction of an automorphism to the algebraically
closed substructure $a(s)$.

A \emph{solution} of $a$ is a functor $\bar{a}: P(n)\rightarrow 
 \mathcal{C}_T$ extending $a$ to the full power set $P(n)$ and satisfying
 the conditions $(i),(ii),(iii)$ (i.e. including the case
 $s=\{1,\dots, n\}$). In particular, in order to find a solution of
 $a$, we need to determine $n$ embeddings 
$$f_i:a(\{1,\ldots,n\}\setminus\{i\})\longrightarrow
 a(\{1,\ldots,n\})=\acl(\{a(i)\mid i \in \{1,\ldots,n\}\}),$$
(for $1\leq i\leq n$) compatible with $a$, that
 is,

$$f_i\circ a_{s,\{1,\ldots,n\}\setminus\{i\}}=f_j\circ
 a_{s,\{1,\ldots,n\}\setminus\{j\}}$$
for every $i,j\in \{1,\ldots,n\}$ and $s\subseteq
\{1,\ldots,n\}\setminus\{i,j\}$. 

The  theory 
 $T$ is said to have \emph{$n$-existence} (over $\acl(\emptyset)$) if every
 $n$-amalgamation problem over $\acl(\emptyset)$ has at least one
 solution. Similarly, we shall
 say that the  theory $T$ has \emph{$n$-uniqueness} (over
 $\acl(\emptyset)$) if every 
 $n$-amalgamation  problem over $\acl(\emptyset)$ has at most one solution up to
 isomorphism (for more details  see~\cite{Hr} and~\cite{Ko}).  
  
It is a well known fact that every simple 
 theory has $2$-existence, by the presence of non-forking
 extensions. Moreover, if the theory is stable, then, by stationarity of
 strong types, $2$-uniqueness holds. Consequentially, also
 $3$-existence holds (for a proof see Lemma~$3.1$ of
 \cite{Hr}). However,  $3$-uniqueness  and $4$-existence can fail 
 for a general stable theory.  Indeed,   in
 \cite{PKM}, the authors thank E. Hrushovski  for supplying an example 
of a stable theory which has neither
$4$-existence nor $3$-uniqueness. The example is the following. Its construction involves a finite cover (for more details about finite covers see \cite{EMI}).
 
\begin{ex}\label{esempio1}Let $\Omega$ be a countable set,
    $[\Omega]^2$  the set of  $2$-subsets of $\Omega$, and
    $C=[\Omega]^2\times \Z/2\Z$. Also let
    $E\subseteq \Omega\times [\Omega]^2$ be the membership relation, and let
    $P$ be the subset of $C^3$ such that
    $((w_1,\delta_1),(w_2,\delta_2),(w_3,\delta_3))$ lies in $P$ if and only if
    there are distinct $c_1,c_2,c_3\in \Omega$ such that
    $w_1=\{c_2,c_3\},w_2=\{c_1,c_3\},w_3=\{c_1,c_2\}$ and 
    $\delta_1+\delta_2+\delta_3=0$. Now let $M$ be the model with the
    $3$-sorted universe $\Omega,[\Omega]^2,C$ and equipped with relations $E,P$
    and projection on the first coordinate $\pi:C\rightarrow
    [\Omega]^2$. Since $M$ is a reduct of $(\Omega,
    \Z/2\Z)^{\textrm{eq}}$, we get that $T=\Th (M)$ is stable. It is
    shown in~\cite{PKM} 
    that $T$ has neither  $4$-existence nor $3$-uniqueness.
 \end{ex}
 
In this paper we generalize this
example. We summarize our main results in the following theorem.
\begin{theorem}\label{summary}
For any 
$n\geq 2$, there exists a stable theory $T_n$ such that $T_n$ has
 $(k+1)$-existence and $k$-uniqueness for any $k\leq n$, but $T_n$ has
  neither $(n+2)$-existence nor $(n+1)$-uniqueness.
\end{theorem}
Also in Proposition~\ref{coincides} we prove that, for $n=2$, the
stable theory $T_2$ given in Theorem \ref{summary} coincides with 
the theory in Example~\ref{esempio1}.\\

All the material we present is expressed in a purely algebraic
terminology. Indeed, the problem of $n$-uniqueness for a theory has also a
natural formulation in terms of permutation groups, as  is shown
in~\cite[Proposition~$3.5$]{Hr}. We adopt this approach here.

In Section~2, we introduce certain permutation modules which  will be
used to construct the automorphism groups of the countable
$\aleph_0$-categorical structures $M_n$ on which is based Theorem~\ref{summary}.  

As  is clear from the definition, the study of amalgamation problems
requires a precise understanding of the algebraic closure in
$T^{\textrm{eq}}$. Since the structures $M_n$ are countable and
$\aleph_0$-categorical, the algebraic closure can
be rephrased with group theoretic terminology: it can be determined
by studying certain closed subgroups of the automorphism group of
$M_n$. This is done in Section~3 and Section~4.

\section{The $\Sym(\Omega)$-submodule structure of
   $\mathbb{F}^{[\Omega]^n}$}\label{sec2} 
   We begin by reviewing some definitions and basic facts about
   permutation groups and permutation modules. 
   
   If $C$ is a set, then the symmetric group $\Sym(C)$  on $C$
   can be considered as a topological group. The
   open sets in this topology are arbitrary unions of cosets of pointwise 
   stabilizers of finite subsets of $C$.  A subgroup $\Gamma$ of
   $\Sym(C)$ is closed if and only if each element of
   $\Sym(C)$ which preserves all the orbits of  $\Gamma$ on  $C^n$, for all
   $n\in \N$, is in $\Gamma$.  It  is well known that  closed subgroups in
   this topology are precisely automorphism groups of first-order
   structures on $C$, see~\cite[Theorem~$5.7$]{Ca} or~\cite{Au}.  
   
Throughout the sequel we denote by $\Fi$ a field, $\Fi_2$ 
the integers modulo $2$, $\Omega$ a countable set and  $[\Omega]^n$ the
set of $n$-subsets of $\Omega$. 
 
 The natural action of the symmetric group $\Sym(\Omega)$ on $[\Omega]^n$ turns
 $\mathbb{F}[\Omega]^n$, the vector 
space over $\Fi$ with basis consisting of the elements of
$[\Omega]^n$, into a $\Sym(\Omega)$-module. We will characterize the
submodules of $\mathbb{F}[\Omega]^n$ in terms of certain
$\Sym(\Omega)$-homomorphisms. The following definition is based on
concepts first introduced in~\cite{Jam}.

\begin{defn}[\cite{Gr}, Def.~$3.4$]\label{beta map}
If $0\leq j\leq n$, then the map $\beta_{n,j}:\mathbb{F}[\Omega]^n\to
\mathbb{F}[\Omega]^j$, given by
$$
\beta_{n,j}(\omega)=\sum_{\omega'\in[\omega]^j}\omega'\qquad
(\textrm{for }\omega\in[\Omega]^n) 
$$
and extended linearly to $\mathbb{F}[\Omega]^n$, is a
$\Sym(\Omega)$-homomorphism (in here we 
denote by $[\omega]^j$  the set of $j$-subsets of $\omega$).  
\end{defn}

It  is shown in
\cite{Gr} (see also \cite{Jam}) that the submodules of
$\mathbb{F}[\Omega]^n$ are completely determined by the  maps $\beta_{n,j}$.
Indeed, it is proved  in \cite[Corollary~$3.17$]{Gr} that every
submodule $U$ of $\mathbb{F}[\Omega]^n$ is an  intersection of  kernels of
$\beta$-maps, i.e. $U=\cap_{j\in S}\ker\beta_{n,j}$ for some subset
$S$ of $\{0,\ldots,n\}$.  

Using the controvariant Pontriagin duality we have that  the dual module of
$\mathbb{F}[\Omega]^n$  is $\Fi^{[\Omega]^n}$, i.e. the set of functions
from $[\Omega]^n$ to 
$\Fi$. We recall that $\Fi^{[\Omega]^n}$ has a natural faithful
action on $[\Omega]^n\times \Fi$ given by
$(w,\delta)^f=(w,f(w)+\delta)$. Hence,  $\Fi^{[\Omega]^n}$, endowed
with the relative topology, 
becomes  a  topological $\Sym(\Omega)$-module and a profinite subgroup
of $\Sym([\Omega]^n\times \Fi)$. Also, given any map
$\beta_{n,j}:\mathbb{F}[\Omega]^n \rightarrow   \mathbb{F}[\Omega]^j 
$, there is a natural  dual  continuous $\Sym(\Omega)$-homomorphism 
$\beta^\ast_{n,j}:\mathbb{F}^{[\Omega]^j}\rightarrow \mathbb{F}^{[\Omega]^n}$
defined by
$$
(\beta^\ast_{n,j}f)(\omega)=\sum_{x\in [\omega]^j}f(x).
$$

Now, the lattice of the closed submodules of  $\mathbb{F}^{[\Omega]^n}$ is the
dual of the lattice of the submodules of $\mathbb{F}[\Omega]^n$.  We
point out that using the algorithm 
described in \cite[Section~$5$]{Gr}, the lattice of the closed submodules of
$\Fi^{[\Omega]^n}$ can be easily computed. Here we record the
following fact that we are frequently going to use. 

\begin{prop}\label{exact}
For $n\geq 1$, $\mathbb{F}=\mathbb{F}_p$ with $p>0$, we have $\Imm \beta_{n,n-1}^\ast=\Ker \beta_{n+1,n}^\ast$.
\end{prop}
\emph{Proof.} The submodule  $\Imm \beta_{n+1,n}$ of
$\mathbb{F}[\Omega]^{n}$ is of the form  $\cap_{j\in S}\Ker
\beta_{n,j}$, for some subset $S$ of
$\{0,\ldots,n\}$. By~\cite[Proposition~$3.19$]{Gr}, we have that $\Imm 
\beta_{n+1,n}\subseteq \Ker \beta_{n,j}$ if and only if $2$ divides
$n+1-j$. Therefore $S=\{j \mid  2 \textrm{ divides } n+1-j\}$.

Also by~\cite[Proposition~$4.1$]{Gr}, we have that if $2$ divides
$n+1-j$, then $\Ker \beta_{n,n-1}\subseteq \Ker \beta_{n,j}$. This yields  
$\Imm \beta_{n+1,n}=\cap_{j\in S}\Ker\beta_{n,j}= \Ker
\beta_{n,n-1}$. In particular, the sequence

\begin{displaymath}
\xymatrix{
\mathbb{F}[\Omega]^{n+1}
\ar@{->}[rr]^{\beta_{n+1,n}}
&&\mathbb{F}[\Omega]^n
\ar@{->}[rr]^{\beta_{n,n-1}}
&&\mathbb{F}[\Omega]^{n-1}}
\end{displaymath}
is exact.

Now the Pontriagin duality is an exact controvariant functor on the
sequences of the 
form $A\to B\to C$. This says that $\Imm \beta_{n,n-1}^\ast=\Ker
\beta_{n+1,n}^\ast$.
\cvd 

\section{Closed submodules of finite index in
  $\mathbb{F}_2^{[\Omega]^n}$ }\label{sec3}

If $A$ is a finite subset of $\Omega$, then we write simply
$\Sym(\Omega\setminus A)$ for the subgroup of $\Sym(\Omega)$ fixing
pointwise $A$. In this section we study the closed
$\Sym(\Omega\setminus A)$-submodules
of $\mathbb{F}_2^{[\Omega]^{n-1}}$ of finite index. We start by
considering the case $A=\emptyset$. 

\begin{lemma}\label{finite index}
If $n\geq 1$, then $\mathbb{F}_2^{[\Omega]^n}$ has no
proper closed $\Sym(\Omega)$-submodule of finite index. 
\end{lemma}

\emph{Proof.}Let $K$ be a closed submodule of
  $\Fi_2^{[\Omega]^n}$ of finite index. Then, 
  $\Fi_2^{[\Omega]^n}/K$ is a finite $\Sym(\Omega)$-module. Since
  $\Sym(\Omega)$ has no proper subgroup of finite index, we get that 
  $\Sym(\Omega)$ centralizes $\Fi_2^{[\Omega]^n}/K$. It follows that
  $f^\sigma -f\in K$, for every $\sigma\in  \Sym(\Omega)$. 

Let $L$ be the annihilator of $K$ in $\mathbb{F}_2[\Omega]^n$,
i.e. $L=\{w\in \mathbb{F}_2[\Omega]^n\mid  g(w)=0 \textrm{ for every
}g\in K\}$. Since $K$ is a closed $\Sym(\Omega)$-submodule,  the set $L$
is a $\Sym(\Omega)$-submodule of 
$\mathbb{F}_2[\Omega]^n$. Now, let $f$ be in
$\mathbb{F}_2^{[\Omega]^n}$, $\sigma$ in $\Sym(\Omega)$ and $w$ in
$L$. We
get $$0=(f^\sigma-f)(w)=f^{\sigma}(w)-f(w)=f(w^{\sigma^{-1}}-w).$$  
This says that $w^{\sigma^{-1}}-w$ is annihilated by every
element of $\mathbb{F}_2^{[\Omega]^n}$. Therefore,
$w^{\sigma^{-1}}-w=0$ and $\sigma$ centralizes  $w$. This shows
that $\Sym(\Omega)$ centralizes $L$. Since $n\geq 1$, the only element
of $\mathbb{F}_2[\Omega]^n$ centralized by $\Sym(\Omega)$ is the zero
vector. Hence $L=0$ and, by the Pontriagin duality,
$K=\mathbb{F}_2^{[\Omega]^n}$. 
\cvd

In the forthcoming analysis  we shall denote  finite subsets of $\Omega$
by capital letters, while  the elements of $[\Omega]^n$ will be generally
denoted by lower cases.   

Now, let $A$ be a finite subset of $\Omega$. To
describe the closed $\Sym(\Omega\setminus A)$-submodules 
of
$\mathbb{F}_2^{[\Omega]^{n-1}}$ of finite
index  we have to
introduce some notation. Let $B$ be a subset of $A$.  We denote  by
$V_{B,A}$ the $\Sym(\Omega\setminus A)$-submodule of 
$\Fi_2^{[\Omega]^{n-1}}$ defined by

\begin{equation}\label{eq:1}
 V_{B,A} = \{f\in \Fi_2^{[\Omega]^{n-1}}  \mid f(w)=0 \, \, \forall \, w\in [\Omega]^{n-1} \textrm{ with } w\cap A\neq B\}
\end{equation}

and we denote by $V_A$ the $\Sym(\Omega\setminus A)$-submodule of
$\Fi_2^{[\Omega]^{n-1}}$ defined by
\begin{eqnarray}\label{eq:2}
V_A&=&\bigoplus_{
\small
B\subseteq A,
|B|< n-1
\normalsize
}V_{B,A}.
\end{eqnarray} 

In the following lemma we describe the elements of $V_A$.

\begin{lemma}\label{eltsVA}Let $A$ be a finite subset of
  $\Omega$. Then 
\begin{equation}\label{eq:3}
V_A=\{f\in\Fi_2^{[\Omega]^{n-1}}\mid f(w)=0\textrm{
    for every }w\in [A]^{n-1}\}.
\end{equation}
\end{lemma}
\emph{Proof.}
We denote by $W$ the vector space on the right hand side of
Equation~$(\ref{eq:3})$. We start by proving that $V_A\subseteq W$. 
Let $B$ be a subset of $A$ with $|B|<n-1$ and  $f$ be in
$V_{B,A}$. Consider $w$ in $[A]^{n-1}$. Since $|B|<n-1$, $|w|=n-1$ and
$w\subseteq A$, we have $w\cap A=w\neq 
B$. By Equation~$(\ref{eq:1})$, we get $f(w)=0$. 
This implies $f\in W$ and so
$V_{B,A}\subseteq W$. Thence, by Equation~$(\ref{eq:2})$, we obtain
$V_A\subseteq W$.

Conversely, we prove that $W\subseteq V_A$. Let $f$ be in $W$. For every subset
$B$ of $A$ with $|B|<n-1$ define
\[
f_B(w)=\left\{
\begin{array}{ccl}
f(w)&&\textrm{if }w\cap A=B,\\
0 &&\textrm{if }w\cap A\neq B.
\end{array}
\right.
\]
Clearly, $f_B\in \Fi_2^{[\Omega]^{n-1}}$ and, by Equation~$(\ref{eq:1})$,
$f_B\in V_{B,A}$. Let $w$ be in $[\Omega]^{n-1}$ with $w\nsubseteq
A$. Since $|w\cap A|<n-1$, we have
\begin{eqnarray*}
\left(\sum_{
\small
B\subseteq A,
|B|< n-1
\normalsize
}f_B\right)(w)&=&\sum_{
\small
B\subseteq A,
|B|< n-1
\normalsize
}f_B(w)=f_{w\cap A}(w)=f(w).
\end{eqnarray*}
Similarly, let $w$ be in $[\Omega]^{n-1}$
with $w\subseteq  A$ (that is, $w\in [A]^{n-1}$). As $f\in W$, we have
$f(w)=0$. Also, by definition of $f_B$, we obtain $f_B(w)=0$. This
shows that $f=\sum_{B\subseteq 
  A,|B|<n-1}f_B$. By Equation~$(\ref{eq:2})$, it follows that
$f\in V_A$.  
\cvd

\begin{lemma}\label{sommadiretta}
Let $A$ be a finite subset of $\Omega$. For each $B\subseteq A$, the
$\Sym(\Omega\setminus 
A)$-modules $V_{B,A}$ are closed submodules of 
$\Fi_2^{[\Omega]^{n-1}}$. Moreover,  
\begin{equation}\label{eq:1111}
\Fi_2^{[\Omega]^{n-1}}=\bigoplus_{
\small
B\subseteq A,
|B|\leq n-1
\normalsize
}V_{B,A}
\end{equation}
and each  $V_{B,A}$ is  $\Sym(\Omega\setminus A)$-isomorphic to
$\Fi_2^{[\Omega\setminus A]^{n-1-|B|}}$. 
\end{lemma}

\emph{Proof.}Since $V_{B,A}$ is an intersection of
  pointwise stabilizers of finite sets of $[\Omega]^{n-1}\times
  \Fi_2$, it is  closed in 
  $\Fi_2^{[\Omega]^{n-1}}$.  It is 
  straightforward to verify the remaining statements.
\cvd

\begin{lemma}\label{primainclusione}
Let $A$ be a finite subset of $\Omega$.
The module $V_A$ has finite index in
$\mathbb{F}_2^{[\Omega]^{n-1}}$. Also, if
$V$ is a closed  $\Sym(\Omega\setminus A)$-submodule 
of $\Fi_2^{[\Omega]^{n-1}}$ of finite index, then $V_A\subseteq V$.
\end{lemma}
\emph{Proof.} By Equations~$(\ref{eq:2})$ and~$(\ref{eq:1111})$,
  we have that $\mathbb{F}_2^{[\Omega]^{n-1}}/V_A$ is isomorphic to
  $\oplus_{|B|=n-1}V_{B,A}$, which has dimension ${|A|\choose
    n-1}$. Therefore $V_A$ has finite index in
  $\mathbb{F}_2^{[\Omega]^{n-1}}$.

Let $V$ be a closed $\Sym(\Omega\setminus A)$-submodule
 of $\mathbb{F}_2^{[\Omega]^{n-1}}$ of finite index. Let $B\subseteq A$ with
$|B|<n-1$. By  Lemma 
  \ref{sommadiretta},  $V_{B,A}$ is
  $\Sym(\Omega\setminus A)$-isomorphic to $\Fi_2^{[\Omega\setminus
      A]^{n-1-|B|}}$. Since $\lbrack V_{B,A}:V_{B,A}\cap V \rbrack=\lbrack
  V_{B,A}+V: V\rbrack$ is finite, we have that $V_{B,A}\cap V$ has
  finite index in 
  $V_{B,A}$. Now, by Lemma~\ref{finite index}, the module   $V_{B,A}$
  does not have 
  any proper closed  $\Sym(\Omega\setminus A)$-submodule of finite
  index. Therefore $V_{B,A}=V_{B,A}\cap V$ and $V_{B,A}\subseteq V$. By
  definition of $V_A$ in Equation~$(\ref{eq:2})$, we get $V_A\subseteq V$.
\cvd
 
In the following lemma we describe the elements of $V_A+\Ker\beta_{n,n-1}^\ast$.

\begin{lemma}\label{char}Let $A$ be a finite subset of $\Omega$. We
  have $V_A+\Ker\beta_{n,n-1}^\ast=\{f\in
  \mathbb{F}_2^{[\Omega]^{n-1}}\mid (\beta_{n,n-1}^\ast f)(w)=0
  \textrm{ for every }w\in [A]^n\}$.
\end{lemma}
\emph{Proof.}If $n=1$, then the equality is clear. So assume $n\geq 2$.

By Lemma~\ref{eltsVA}, the elements of $V_A$ are the
  functions $f\in 
  \mathbb{F}_2^{[\Omega]^{n-1}}$ vanishing on each element of $
         [A]^{n-1}$. Now, if $f_1\in V_A$, $f_2\in \Ker\beta_{n,n-1}^\ast$ and
         $w\in [A]^n$,
         then $$(\beta_{n,n-1}^\ast (f_1+f_2))(w)=(\beta_{n,n-1}^\ast
         f_1)(w)=\sum_{w'\in[w]^{n-1}}f_1(w')=0. $$

Therefore, it remains to prove that if $f\in
\mathbb{F}_2^{[\Omega]^{n-1}}$ and $(\beta_{n,n-1}^\ast f)(w)=0$ for
every $w\in [A]^{n}$, then $f\in V_A+\Ker\beta_{n,n-1}^\ast $. Let $a$
be a fixed element of $A$ and let
$g\in \mathbb{F}_2^{[\Omega]^{n-2}}$ be the function defined by   
\[
g(\omega)=\left\{
\begin{array}{cll}
f(\omega\cup \{a\})&&\textrm{if }\omega\subseteq A \textrm{ and
}a\notin \omega,\\ 
0&&\textrm{otherwise }.
\end{array}
\right.
\]
Set $f_2=\beta_{n-1,n-2}^\ast g$. By Proposition~\ref{exact}, we have that
$f_2\in \Imm\beta_{n-1,n-2}^\ast=\Ker\beta_{n,n-1}^\ast$. Set
$f_1=f-f_2$. We claim that $f_1$
lies in $V_A$, from which the lemma follows. By Lemma~\ref{eltsVA}, it
suffices to prove that  
$f_1(w')=0$ for every $w'\in [A]^{n-1}$. Let $w'$ be in
$[A]^{n-1}$. Assume $a\in w'$. By the definition of $g$, we have
$$f_2(w')=(\beta_{n-1,n-2}^\ast g)(w')=\sum_{\omega\in
  [w']^{n-2}}g(\omega)=g(w'\setminus \{a\})=f(w')$$
and $f_1(w')=0$. Now assume $a\notin w'$. By the definition of $g$ and by
the hypothesis on $f$, we have
\begin{eqnarray*}
f_2(w')&=&(\beta_{n-1,n-2}^\ast g)(w')=\sum_{\omega\in
  [w']^{n-2}}g(\omega)=\sum_{\omega\in [w']^{n-2}}f(\omega\cup
\{a\})\\
&=&\sum_{x\in
  [w'\cup\{a\}]^{n-1}}\!\!\!\!\!\!f(x)+f(w')=(\beta_{n,n-1}^\ast
f)(w'\cup\{a\})+f(w')=f(w'), 
\end{eqnarray*}
and $f_1(w')=0$.
\cvd

\begin{defn}\label{defn}
We write  $W_A$ for $\beta_{n,n-1}^\ast(V_A)$, with $V_A$ as in
Equation~$(\ref{eq:2})$. 
\end{defn}
Now, using the previous
lemmas we describe the 
closed $\Sym(\Omega\setminus  A)$-submodules  of
$\Imm\beta_{n,n-1}^\ast$ of finite index.  

\begin{prop}\label{uniqueimage}Let $A$ be a finite subset of
  $\Omega$. The module $W_A$ is the unique minimal
 closed $\Sym(\Omega\setminus A)$-submodule of
 $\Imm\beta_{n,n-1}^\ast$ of finite index. Furthermore, 
  $W_A=\{g\in \Imm\beta_{n,n-1}^\ast\mid g(w)=0 \textrm{ for every } w\in
  [A]^n\}$.
\end{prop}
\emph{Proof.}Let $W$ be a closed $\Sym(\Omega\setminus A)$-submodule
  of $\Imm\beta_{n,n-1}^\ast$ of finite index. By the first
  isomorphism theorem $W$ is the image via $\beta_{n,n-1}^\ast$ of some closed
  $\Sym(\Omega\setminus A)$-submodule $V$ of
  $\mathbb{F}_2^{[\Omega]^{n-1}}$ of finite index. Now, by
  Lemma~\ref{primainclusione}, we get 
  $V_A\subseteq V$. So
  $\beta_{n,n-1}^\ast(V_A)\subseteq \beta_{n,n-1}^\ast(V)= 
  W$. Hence, $W_A=\beta_{n,n-1}^\ast(V_A)$ is the unique minimal
  closed $\Sym(\Omega\setminus A)$-submodule of
  $\Imm\beta_{n,n-1}^\ast$ of finite index. 

Now, from Lemma~\ref{char} the rest of the proposition is immediate. 
\cvd

\section{The infinite family of examples}\label{sec1}

Before introducing our examples, we need to set some
auxiliary notation. 

\begin{defn}\label{auxiliary} Let $M$ be a structure and $A,B$ subsets
  of $M$. We denote 
  by $\overline{\Aut(A/B)}$ the subgroup of $\Aut(M)$
  fixing setwise $A$ and fixing pointwise $B$. The setwise stabilizer of $A$ in $\Aut(M)$ will be denoted by $\Aut(M)_{\{A\}}$, while the permutation group
  induced by $\overline{\Aut(A/B)}$ on $A$ will be denoted by $\Aut(A/B)$. 
\end{defn}

Let $n\geq 2$ be an integer and $\Omega$ be a countable set. 

\begin{defn}\label{OURdefinition}We
consider $M_n$  
 the multisorted structure with 
sorts $\Omega$, $[\Omega]^n$ and $[\Omega]^n\times \mathbb{F}_2$ and with 
automorphism group $\Imm \beta^\ast_{n,n-1}\rtimes \Sym(\Omega)$. Note 
that this is well-defined as $\Imm\beta_{n,n-1}^\ast$ is a closed
submodule of $\mathbb{F}_2^{[\Omega]^n}$.

Moreover, the theory $T_n=\Th (M_n)$ is stable (see Section \ref{coincidence}).  
\end{defn}

In the next paragraph we introduce some notation that would
be useful to describe the  algebraically closed sets of $M_n$. 

Denote by
$\pi:[\Omega]^n\times \mathbb{F}_2\to [\Omega]^n$ the projection on
the first coordinate. Given  $A$ a finite subset of  $M_n$, we have
that $A$ is of the form $A_1\cup A_2\cup A_3$, where $A_1$ belongs  to
the sort $\Omega$, $A_2$ belongs to the sort $[\Omega]^n$ and $A_3$
belongs to the sort $[\Omega]^n\times \mathbb{F}_2$. Consider
$\tilde{A_2}\subseteq 
\Omega$ the union of the elements in $A_2$ and $\tilde{A_3}\subseteq
\Omega$ the union of the elements in $\pi(A_3)$. We define
the \emph{support} of $A$, written $\supp(A)$, to be the subset $A_1\cup
\tilde{A_2}\cup \tilde{A_3}$ of $\Omega$. Finally, we define $\cl(A)$ to be the subset of $M_n$ 
$$
\cl(A):=\supp(A)\cup 
 [\supp(A)]^n\cup  ([\supp(A)]^n\times \mathbb{F}_2)
 $$

\smallskip

In the rest of this section we  describe the algebraically closed sets in the
structure $M_n$. Here we  consider structures \emph{up to
  interdefinability}, which allows us to identify an
$\aleph_0$-categorical structure with its automorphism group. So we
identify two substructures $A_1,A_2$ of a structure $M$, if
$\Aut(A_1)=\Aut(A_2)$. If $M$ is an $\aleph_0$-categorical structure and
$A\subset M$, we denote the algebraic closure $\acl^{\rm{eq}}(A)$ of
$A$ simply by $\acl(A)$, i.e. 
the union of the finite  $\Aut(M/A)$-invariant sets of
$M^{\rm{eq}}$. We recall that definable subsets of $\acl(A)$
correspond, up to interdefinability,  to closed subgroups of
$\Aut(M/A)$ of finite
index, see~\cite[Section~$4.1$]{Wi} or Theorem~$4.1$ in the article
\emph{``The structure of totally   categorical structures''} by
W. Hodges~\cite[page~$116$]{Au}.  

Similarly, if $A\subset M$, we denote the definable closure $\dcl^{\rm{eq}}(A)$ of
$A$ simply by $\dcl(A)$, i.e. 
the set of the points of $M^{\rm{eq}}$ fixed by $\Aut(M/A)$. 

\begin{lemma}\label{conclusione1}
Let $A$ be a finite set of $M_n$. Then 
$$\Aut(M_n/\cl(A))=W_{\supp(A)}\rtimes \Sym(\Omega\setminus \supp(A))$$
 (where $W_{\supp(A)}$ is the closed $\Sym(\Omega\setminus\supp(A))$-submodule of
$\Imm\beta_{n,n-1}^\ast$ in Definition~\ref{defn}). Moreover, $\Aut(M_n/\cl(A))$ is the unique
minimal closed subgroup of finite index of $\Aut(M_n/A)$.
\end{lemma}
\emph{Proof.}
Set $\Gamma=\Aut(M_n/\cl(A))$. We first prove that 
$\Gamma=W_{\supp(A)}\rtimes \Sym(\Omega\setminus \supp(A))$. By
definition of the multisorted structure $M_n$, we have $\Aut
M_n=\Imm\beta_{n,n-1}^*\rtimes 
\Sym(\Omega)$. Therefore, an element of $\Gamma$ is
an ordered pair of the form $g\sigma$, where $g\in
\Imm\beta_{n,n-1}^*$ and $\sigma\in 
\Sym(\Omega)$. The action of $g\sigma$ on the elements  belonging
to the sorts $\Omega$ and $[\Omega]^n$ is given by the permutation
$\sigma$. Also, the action of $g\sigma$ on the element $(w,x)$ belonging to
the sort $[\Omega]^n\times \Fi_2$ is given
by $$(w,x)^{g\sigma}=(w^\sigma,x+g(w)).$$ 
This implies that  the automorphism $g\sigma$ fixes the elements in 
$\supp(A)$ and in $[\supp(A)]^n$ (in the sorts $\Omega$ and
$[\Omega]^n$) if and only if 
$\sigma\in\Sym(\Omega\setminus\supp(A))$.  Also, the automorphism
$g\sigma$ fixes the
elements in $[\supp(A)]^n\times \Fi_2$ (in the sort
$[\Omega]^n\times \Fi_2$) if and only if $g(w)=0$ for
every $w\in [\supp(A)]^n$. Hence, by the description of the elements of
$W_{\supp(A)}$ in Proposition~\ref{uniqueimage}, we have $g\sigma\in
\Gamma$ if and only if $g\sigma\in W_{\supp(A)}\rtimes
\Sym(\Omega\setminus \supp(A))$. 

 We claim that $\Gamma$ is the unique
minimal closed subgroup of $\Aut(M_n/A)$ of finite index. Note that $\Gamma$ is a closed subgroup of
$\Aut(M_n/A)$ of finite index.

Now, let $H$ be
a closed subgroup of 
$\Aut(M_n/A)$ 
of finite index. Up to
replacing $H$ with $H\cap \Gamma$, we may assume that $H\subseteq
\Gamma$. Let
$\mu:\Gamma\to\Sym(\Omega\setminus \supp(A))$ be the 
natural projection. Since $\mu$ is a surjective continuous closed map
and $\Sym(\Omega\setminus \supp(A))$ has no proper subgroup of finite
index, we get that $\mu(H)=\Sym(\Omega\setminus \supp(A))$. This
yields that $H\cap W_{\supp(A)}$ is a closed
$\Sym(\Omega\setminus\supp(A))$-submodule of $W_{\supp(A)}$ of finite
index. Now Proposition~\ref{uniqueimage} shows that $H\cap
W_{\supp(A)}=W_{\supp(A)}$. So $W_{\supp(A)}\subseteq H$ and $H=\Gamma$. 
\cvd
In the following we denote by $\acl_{M_n}$ the $\acl$ in $M_n$.
\begin{prop}\label{acl1} Let $A$ be a finite set of $M_n$. Then
$\acl_{M_n}(A)=\cl(A)$.  
\end{prop}
\emph{Proof.} Let $\overline{b}$ be an $m$-tuple in $M_n$ and $A$ be a finite set of $M_n$. We first claim that $\Aut(M_n/\overline{b})\geq \Aut(M_n/\cl(A))$ if and only if the underlying set of $\overline{b}$ is conteined in $\cl(A)$ . One direction is obvious. Suppose that $\Aut(M_n/\overline{b})\geq \Aut(M_n/\cl(A))$ for some finite $A\subset M_n$. Then by Lemma \ref{conclusione1} we have that $\Aut(M_n/\cl(\cl(A),\overline{b})$ is a closed subgroup of finite index in  $\Aut(M_n/\cl(A),\overline{b})=\Aut(M_n/\cl(A))$. Hence  $\Aut(M_n/\cl(\cl(A),\overline{b})$ is a closed subgroup of finite index in $\Aut(M_n/A)$. By uniqueness of the minimal closed subgroup of finite index of $\Aut(M_n/A)$ we get that $W_{\supp(A)}\rtimes \Sym(\Omega\setminus \supp(A))$ is equal to $W_{\supp(\cl(A),\overline{b})}\rtimes \Sym(\Omega\setminus \supp(\cl(A), \overline{b}))$ and, since $\supp(\cl(A), \overline{b})=\supp(A,\overline{b})$, this is possible if and only if $\supp(\overline{b})\subseteq \supp(A)$, which proves the claim.

By definition, $\acl_{M_n}(A)$ is the union of the finite orbits on $M_n$ of $\Aut(M_n/A)$. Let $c\in \acl_{M_n}(A)$. Then $\Aut(M_n/A,c)$ is a closed subgroup of  finite index in $\Aut(M_n/A)$. Hence, by Lemma \ref{conclusione1}, $\Aut(M_n/A,c)\geq \Aut(M_n/\cl(A)$. By the above argument we have that $c\in \cl(A)$. 

Let $c\in \cl(A)$, then $\Aut(M_n/A)\geq\Aut(M_n/A,c)\geq \Aut(M_n/\cl(A))$. Hence the index of $\Aut(M_n/A,c)$ in $\Aut(M_n/A)$ is finite. 
\cvd
 Let  $c^{\rm{eq}}\in M_n^{\rm{eq}}$. Then $c^{\rm{eq}}$  is a $0$-definable equivalence class of a  tuple $b$ of elements in $M_n$.  We denote by $\mathcal{s}(c^{\rm{eq}})$ the union of elements in $M_n$ of  $c^{\rm{eq}}$. Similarly if $A\subseteq M_n^{\rm{eq}}$, we denote by $\mathcal{s}(A)$ the set of elements in $M_n$ $\bigcup_{c^{\rm{eq}}\in A}\mathcal{s}(c^{\rm{eq}})$.
\begin{prop}\label{acl2} Let $A$ be a finite set of $M_n$. Then $\mathcal{s}(\acl(A))=\cl(A)$. In particular 
 $\acl(\emptyset)=\emptyset$.
\end{prop}
\emph{Proof.} Fix an enumeration $\overline{b}$ of $\acl_{M_n}(A)$ and set  $\Gamma=\Aut(M_n/\acl_{M_n}(A))$.  Consider the trivial relation $R=\{(b^\alpha,b^\alpha):\alpha\in \Aut(M_n)\}$. Since $R$ is an $\Aut(M_n)$-orbit, $R$ is a $0$-definable equivalence relation in $M_n$. Consider the $R$-equivalence class of $\overline{b}$. The pointwise stabilizer of $\overline{b}$ in $\Aut(M_n)$ is $\Gamma$ which, by Lemma \ref{conclusione1} and Proposition \ref{acl1}, has finite index in $\Aut(M_n/A)$ and so $\overline{b}\in \acl(A)$.

Let $c^{\rm{eq}}\in \acl(A)$, then $\Aut(M_n/A,c^{\rm{eq}})$ is a closed subgroup of finite index of $\Aut(M_n/A)$. By Lemma \ref{conclusione1} $\Aut(M_n/A,c^{\rm{eq}})$ contains $\Gamma$. Being  $\Aut(M_n/A,c^{\rm{eq}})$ also open in $\Aut(M_n/A)$ there exists a finite tuple $\overline{b}$ in $M_n$ such that $\Aut(M_n/A,c^{\rm{eq}})$ contains the basic open subgroup $\Aut(M_n/A, \overline{b})$. Moreover $c^{\rm{eq}}=\overline{b}^{\Aut(M_n/A,c^{\rm{eq}})}$. By $\aleph_0$-categoricity  the index of $\Aut(M_n/A, \overline{b})$ in $\Aut(M_n/A,c^{\rm{eq}})$ is finite. Then, the index  of $\Aut(M_n/A, \overline{b})$ in $\Aut(M_n/A)$ is finite and so $\Gamma\leq \Aut(M_n/A, \overline{b})$. Hence by the same argument used in Proposition \ref{acl1}, we get that the underlying set in $M_n$ of $\overline{b}$ is contained in $\cl(A)=\acl_{M_n}(A)$. From the fact that $\Aut(M_n/A,c^{\rm{eq}})\leq \Aut(M_n/A)$ and $\overline{b}\in \acl_{M_n}(A)$ it follows immediately that also the underlying set of the $\Aut(M_n/A,c^{\rm{eq}})$-orbit $\overline{b}^{\Aut(M_n/A,c^{\rm{eq}})}$ is contained in $\acl_{M_n}(A)$.
\cvd

\begin{cor}\label{acl7}
Let $A$ be a finite set of $M_n$. Then,
$$\Aut(M_n)_{\{\acl_{M_n}(A)\}}=\Aut(M_n)_{\{\acl(A)\}}.$$
\end{cor}
\emph{Proof.} From Proposition \ref{acl1} and Proposition \ref{acl2} it follows that $\Aut(M_n)_{\{\acl(A)\}}\leq \Aut(M_n)_{\{\acl_{M_n}(A)\}}$. Now, let $g\in  \Aut(M_n)_{\{\acl_{M_n}(A)\}}$. Note that $\acl_{M_n}(A^g)=\acl_{M_n}(A)$. Consequently, $\acl(A^g)=\acl(A)$.  If  $c^{\rm{eq}}\in \acl(A)$, then the index of $\Aut(M_n/A, c^{\rm{eq}})$ in $\Aut(M_n/A)$ is finite. Therefore, $\Aut(M_n/A^g, (c^{\rm{eq}})^g)=g^{-1}\Aut(M_n/A, c^{\rm{eq}})g$ has finite index in $\Aut(M_n/A^g)=g^{-1}\Aut(M_n/A)g$, which implies that $(c^{\rm{eq}})^g\in \acl(A^g)=\acl(A)$.
\cvd

\begin{prop}\label{acl3}
Let $A$ be a finite subset of $M_n$. Then,
$
\dcl(\acl_{M_n}(A))=\acl(A).$
\end{prop}
\emph{Proof.} Let $c^{\rm{eq}}\in \acl(A)$, i.e. the stabilizer of  $c^{\rm{eq}}$ in $\Aut(M_n/A)$ has finite index in $\Aut(M_n/A)$. We need to show that the stabilizer of  $c^{\rm{eq}}$ in $\Aut(M_n/\acl_{M_n}(A))$ is equal to $\Aut(M_n/\acl_{M_n}(A))$. We have the following disequality:
$$
|\Aut(M_n/\acl_{M_n}(A)):\Aut(M_n/\acl_{M_n}(A), c^{\rm{eq}})|\leq |\Aut(M_n/A):\Aut(M_n/A, c^{\rm{eq}})|
$$
Then $|\Aut(M_n/A):\Aut(M_n/\acl_{M_n}(A), c^{\rm{eq}})|$ is finite. By Lemma \ref{conclusione1} and Proposition \ref{acl1} it follows that $\Aut(M_n/\acl_{M_n}(A), c^{\rm{eq}})$, is equal to $\Aut(M_n/\acl_{M_n}(A))$, i.e. $c^{\rm{eq}}\in \dcl(\acl_{M_n}(A))$.

Let $c^{\rm{eq}}\in \dcl(\acl_{M_n}(A))$.  We need to show that $\Aut(M_n/A, c^{\rm{eq}})$, has finite index in $\Aut(M_n/A)$. We have that 
\begin{equation}\label{indici}
\begin{array}{c}
|\Aut(M_n/A):\Aut(M_n/\cl(A)), c^{\rm{eq}})|=\\
|\Aut(M_n/A):\Aut(M_n/A, c^{\rm{eq}})||\Aut(M_n/A, c^{\rm{eq}}):\Aut(M_n/\cl(A), c^{\rm{eq}})|
\end{array}
\end{equation}
Since $c^{\rm{eq}}\in \dcl(\acl_{M_n}(A))$ we have that $\Aut(M_n/\acl_{M_n}(A), c^{\rm{eq}})=\Aut(M_n/\acl_{M_n}(A))$. Lemma \ref{conclusione1} and the equality (\ref{indici}) imply that $|\Aut(M_n/A):\Aut(M_n/A, c^{\rm{eq}})|$ is finite. This proves that $c^{\rm{eq}}\in \acl(A)$ and the proof is complete.
\cvd

\begin{cor}\label{acl5}
Let $A$ be a finite subset of $M_n$. Then
$$\Aut(M_n/\acl_{M_n}(A))=\Aut(M_n/\acl(A)).$$
\end{cor}
\emph{Proof.}
Let $g\in \Aut(M_n/\acl_{M_n}(A))$ and $c^{\rm{eq}}\in \acl(A)$. Proposition \ref{acl3} yields that $(c^{\rm{eq}})^g=c^{\rm{eq}}$, which means that $g\in \Aut(M_n/\acl(A))$. 
It remains to prove that $\Aut(M_n/\acl(A))\leq \Aut(M_n/\acl_{M_n}(A))$. Consider the trivial relation $R$ given by $R=\{ (b,b) : b\in M_n\}$. This is a 0-definable relation. Let $a\in \acl_{M_n}(A)$.  Then  $\{a\}\in M_n^{\rm{eq}}$ and $\Aut(M_n/A, \{a\})=\Aut(M_n/A, a)$ is a closed subgroup of finite index in $\Aut(M_n/A)$. Hence, we can consider that $\acl_{M_n}(A)\subseteq \acl(A)$ and the thesis follows at once. 
\cvd
\begin{rem}\label{acl4}
 Proposition \ref{acl1} yields that if $A$ is a finite set of $M_n$, then $\acl_{M_n}(A)=\acl_{M_n}(\supp(A))$. Therefore, from Proposition \ref{acl3} it follows  that $\acl(A)=\acl(\supp(A))$.
\end{rem}
\begin{prop}\label{acl6}
Let $A_1, \dots, A_n$ be finite subsets in the sort $\Omega$. Then
$$
\acl(\acl(A_1),\dots,\acl(A_n))=\acl(\bigcup_{i=1}^n A_i).
$$
\end{prop}
\emph{Proof.}
Obviously, $\acl(\bigcup_{k=1}^n A_k)\subseteq \acl(\acl(A_1),\dots,\acl(A_n))$.

Let $c^{\rm{eq}}\in \acl(\acl(A_1),\dots,\acl(A_n))$ and set $G=\Aut(M_n/\acl(A_1),\dots,\acl(A_n))$. Then, the pointwise stabilizer $G_{c^{\rm{eq}}}$ has finite index in $G$. By Corollary \ref{acl5} we have that 
$$
G=\bigcap_{i=1}^n W_{A_i}\rtimes \Sym(\Omega \setminus A_i).
$$
Moreover, $G\geq W_{\bigcup_{i=1}^n A_i}\rtimes \Sym(\Omega \setminus \bigcup_{i=1}^n A_i)$
and $G$ is a closed subgroup in $\Aut(M_n/\bigcup_{i=1}^n A_i)$. So, $G$ is a closed subgroup of finite index in $\Aut(M_n/\bigcup_{i=1}^n A_i)$ which implies that also $G_{c^{\rm{eq}}}$ is of finite index in $\Aut(M_n/\bigcup_{i=1}^n A_i)$. Now, $G_{c^{\rm{eq}}}=G\cap \Aut(M_n/\bigcup_{i=1}^n A_i,c^{\rm{eq}})$ and 
$$
\begin{array}{c}
|\Aut(M_n/\bigcup_{i=1}^n A_i):\Aut(M_n/\bigcup_{i=1}^n A_i, c^{\rm{eq}})|=\\
|\Aut(M_n/\bigcup_{i=1}^n A_i):G_{c^{\rm{eq}}}|/|\Aut(M_n/\bigcup_{i=1}^n A_i, c^{\rm{eq}}):G_{c^{\rm{eq}}}|,
\end{array}
$$
i.e. $c^{\rm{eq}}\in \acl(\bigcup_{i=1}^n A_i)$.
\cvd

\section{$k$-existence and $k$-uniqueness for $M_n$}
In this section we prove Theorem \ref{summary}. Note that, up to
  renaming the elements of $\Omega$, we may assume that
  $\Omega=\mathbb{N}$. In the sequel we denote by $[k]$
   the subset $\{1,\ldots,k\}$ of $\N$. Also, given $i\in [k]$, we
   denote by
  $[k]-i$ the set $\{1,\ldots,k\}\setminus \{i\}$.  Finally, we denote
   the theory $\Th
   (M_n)$ by $T_n$.

We start by studying $k$-uniqueness in $T_n$. We first single out
the following technical lemma which would be used in Proposition~\ref{unik}.
\begin{lemma}\label{lemma:aux}
Let $k$ and $n$ be integers, with $k< n$, and 
$A_1,\ldots,A_{k}$ be subsets of $\Omega$. Then 
$$(\dag)\qquad\bigcap_{i=1}^{k}\left(V_{A_i}+\Ker\beta_{n,n-1}^\ast\right)=
\left(\bigcap_{i=1}^{k}V_{A_i}\right)+\Ker\beta_{n,n-1}^{\ast}.$$  
\end{lemma} 
\emph{Proof.}
We denote the left-hand-side of $(\dag)$ by $V_{1,k}$ and the
right-hand-side of $(\dag)$ by $V_{2,k}$ (where the label $k$ is used
in order to remember
the number of intersections). 

We argue by induction on $k$. Note that if $k=0$ or $k=1$, then there
is nothing to prove.  Assume
$(\dag)$ holds for $k$ intersections (where $k\geq 1$)
and that $k+1<n$. In particular, we point out that $n>2$. We prove
that $(\dag)$ holds for $k+1$ intersections.  
Clearly, $V_{2,k+1}\subseteq V_{1,k+1}$. Let $g$ be in $V_{1,k+1}$. We
need to show that $g\in V_{2,k+1}$. By induction hypothesis (on the
sets $A_1,\ldots,A_k$), 
we
have 
\begin{equation}\label{eq:p1}
V_{1,k+1}=
\left(\left(\bigcap_{i=1}^{k}V_{A_i}\right)+\Ker\beta_{n,n-1}^\ast\right)\cap
(V_{A_{k+1}}+\Ker\beta_{n,n-1}^\ast).
\end{equation}
By Equation~(\ref{eq:p1}) and Proposition~\ref{exact}, we have
\begin{equation}\label{eq:p2}
g=g_1+\beta_{n-1,n-2}^\ast h_1=g_2+\beta_{n-1,n-2}^\ast h_2,
\end{equation}
where 
$g_1\in\cap_{i=1}^kV_{A_i}$, $g_2\in V_{A_{k+1}}$ and $h_1,h_2\in
\Fi_2^{[\Omega]^{n-2}}$. We claim that (up to replacing $h_1$ by
$h_1+l$, where $l\in \Ker\beta_{n-1,n-2}^\ast$), we may assume that
$h_1-h_2\in \cap_{i=1}^k V_{A_i\cap A_{k+1}}$.

Let $w$ be an $(n-1)$-subset of $\Omega$ contained in $A_i\cap A_{k+1}$
for some $i=1,\ldots,k$. Since $g_1\in V_{A_i}$ and $g_2\in
V_{A_{k+1}}$, we see that $g_1(w)=g_2(w)=0$. So, from
Equation~(\ref{eq:p2}) we obtain  
$$g(w)=(\beta_{n-1,n-2}^\ast h_1)(w)=(\beta_{n-1,n-2}^\ast h_2)(w),$$ that is, 
$(\beta_{n-1,n-2}^\ast(h_1-h_2))(w)=0$. As $w$ is an arbitrary $(n-1)$-subset of
$A_i\cap A_{k+1}$, Lemma~\ref{char}
 yields $h_1-h_2\in V_{A_i\cap A_{k+1}}+\Ker\beta_{n-1,n-2}^\ast$. As
 $i$ is an arbitrary element in $\{1,\ldots,k\}$, we get 
$$h_1-h_2\in \bigcap_{i=1}^k(V_{A_i\cap
   A_{k+1}}+\Ker\beta_{n-1,n-2}^\ast).$$ 
Since $k+1<n$, we have $k<n-1$ and so we may now apply our
 inductive hypothesis on the sets $A_1\cap A_{k+1},\ldots,A_{k}\cap
 A_{k+1}$. We have
\begin{equation}\label{eq:p3}
h_1-h_2\in \left(\bigcap_{i=1}^kV_{A_i\cap
   A_{k+1}}\right)+\Ker\beta_{n-1,n-2}^\ast.
\end{equation}
From Equation~(\ref{eq:p3}), we get
$h_1-h_2=h+l$, where $h\in \cap_{i=1}^kV_{A_i\cap
A_{k+1}}$ and $l\in \Ker\beta_{n-1,n-2}^\ast$. Set $h_1'=h_1+l$. We have
$$h_1'-h_2=h_1+l-h_2=h\in \cap_{i=1}^kV_{A_i\cap A_{k+1}}$$
and our claim is proved.

Let $t$ be the element of $\Fi_2^{[\Omega]^{n-2}}$ defined by
\[
t(w)=\left\{
\begin{array}{ccl}
h_1(w)&&\textrm{if }w\subseteq A_i\textrm{ for some }i=1,\ldots,k,\\
h_2(w)&&\textrm{if }w\subseteq A_{k+1},\\
0&&\textrm{otherwise}.
\end{array}
\right.
\] 
Note that the function $t$ is well-defined. Indeed, recall that  $n>2$
and note that if $w$ is an
$(n-2)$-subset of $\Omega$ with $w\subseteq A_i\cap A_{k+1}$ (for some
$i=1,\ldots,k$), then $h_1(w)=h_2(w)$ as $h_1-h_2\in V_{A_{i}\cap
  A_{k+1}}$. 

We claim that $g+\beta_{n-1,n-2}^\ast t \in
\cap_{i=1}^{k+1}V_{A_i}$. We have to show that $g+\beta_{n-1,n-2}^\ast
t$ vanishes in $[A_i]^{n-1}$, for $i=1,\ldots,k+1$.   Let
$w$ be an $(n-1)$-subset of $\Omega$ with $w\subseteq A_i$, for some
$i=1,\ldots,k+1$. If $i\leq k$, then we have

$$(g+\beta_{n-1,n-2}^\ast t)(w)=(g_1(w)+\beta_{n-1,n-2}^\ast
h_1(w))+\beta_{n-1,n-2}h_1(w)=0,$$ 
where in the first equality we used Equation~(\ref{eq:p2}) and the fact
that $t$ and $h_1$ coincide in $[A_i]^{n-2}$, and in the second
equality we used that $g_1\in V_{A_i}$. 
Similarly, if $i=k+1$, then  
$$(g+\beta_{n-1,n-2}^\ast t)(w)=(g_2(w)+\beta_{n-1,n-2}^\ast
h_2(w))+\beta_{n-1,n-2}h_2(w)=0,$$ 
where in the first equality we used Equation~(\ref{eq:p2}) and the fact
that $t$ and $h_2$ coincide in $[A_{k+1}]^{n-2}$, and in the second
equality we used that $g_2\in V_{A_{k+1}}$.
 
Finally, as $\beta_{n-1,n-2}^\ast t\in\Ker\beta_{n,n-1}^\ast$, we get that
$g\in V_{2,k+1}$.
\cvd

\begin{prop}\label{unik}
The theory $T_n$ has $k$-uniqueness for every  $k\leq n$. 
\end{prop}
\emph{Proof.}Let $k$ be an integer with $k\leq n$ and
  $a:P(k)^-\to \mathcal{C}_{T_n}$ be a $k$-amalgamation problem. We need
  to show
  that $a$ has at most one solution up to
  isomorphism. Since every stable theory has $1$- and $2$-uniqueness, we
  may assume that $k\geq 3$.  Set
  $\Gamma_1=\Aut(a([k-1])/\cup_{i=1}^{k-1}a([k]-i))$ 
and $\Gamma_2=\Aut(a([k-1])/
\cup_{i=1}^{k-1}a([k-1]-i))$. By~\cite[Proposition~$3.5$]{Hr}, it is
enough to prove 
  that 
\begin{equation}
\Gamma_1=\Gamma_2, 
\end{equation}
i.e. $\overline{\Gamma_1},\overline{\Gamma_2}$ give rise to the same action on
$a([k-1])$ (see Definition~\ref{auxiliary}).   

By Remark \ref{acl4}, the algebraically closed
sets of finite subsets of $M_n$ are of the form $\acl(A)$, for some finite subset $A$ of the sort
$\Omega$. By Corollary \ref{acl7} the setwise stabilizer of $\acl(A)$ in
$\Aut(M_n)$ is simply $(\Sym(\Omega\setminus A)\times\Sym(A))\ltimes
\Imm\beta_{n,n-1}^\ast$. Using
Corollary~\ref{acl5}, we get that the
pointwise stabilizer of 
$\acl(A)$ in $\Aut(M_n)$ is $\Sym(\Omega\setminus A)\ltimes W_A$. 

Let $a(i)=\acl(B_i)$, where $B_i$ are finite subsets of $M_n$ for $1\leq i\leq k$. Set $A_i=\supp(B_i)$, for $1\leq i\leq k$, and
$A=\cup_{i=1}^{k-1}A_i$. Note that by definition of amalgamation problem and by Proposition \ref{acl6}, we
have $a([k-1])=\acl(A)$. Therefore, by the previous paragraph, as
$k\geq 3$, we
get that $\overline{\Gamma_1}$ is equal to
$$
((\Sym(\Omega\setminus A)\times \Sym(A))\ltimes
\Imm\beta_{n,n-1}^\ast)\cap\bigcap_{i=1}^{k-1}(\Sym(\Omega\setminus
((A\cup A_k)\setminus A_i))\ltimes W_{(A\cup A_k)\setminus A_i})
$$
i.e.  
\begin{equation}\label{eq:G1}
\overline{\Gamma_1}=\Sym(\Omega\setminus (A\cup A_k))\ltimes
\bigcap_{i=1}^{k-1}W_{(A\cup A_k)\setminus A_i} 
\end{equation}
and $\overline{\Gamma_2}$ is equal to
$$
((\Sym(\Omega\setminus A)\times \Sym(A))\ltimes
\Imm\beta_{n,n-1}^\ast)\cap\bigcap_{i=1}^{k-1}(\Sym(\Omega\setminus
(A\setminus A_i))\ltimes W_{A\setminus A_i})
$$
i.e.
\begin{equation}\label{eq:G2}
\overline{\Gamma_2}=\Sym(\Omega\setminus A)\ltimes
\bigcap_{i=1}^{k-1}W_{A\setminus A_i}. 
\end{equation}
As $\Sym(\Omega\setminus(A\cup A_k))$ and $\Sym(\Omega\setminus A)$
act trivially on the elements of $\acl(A)$, by
Equations~$(\ref{eq:G1})$ and~$(\ref{eq:G2})$, in order to prove that
$\Gamma_1=\Gamma_2$ it suffices to show that
\begin{equation*}
W_1=\bigcap_{i=1}^{k-1}W_{(A\cup A_k)\setminus A_i}
\quad\textrm{and}\quad W_2=\bigcap_{i=1}^{k-1}W_{A\setminus A_i}
\end{equation*}
induce the same action on $\acl(A)$. Also, $W_1$ and $W_2$ act
trivially on the elements belonging to the sorts $\Omega$ and
$[\Omega]^n$ of $M_n$. Thus, it suffices to study the action of $W_1$
and $W_2$ on the elements of $\acl(A)$ belonging to the sort
$[\Omega]^n\times \Fi_2$, that is, on
$[A]^n$.  Clearly,
$W_1\subseteq W_2$. Therefore, it remains to show that for every element
$f$ of $W_2$ there exists an element $\overline{f}$ of $W_1$ such that
$f$ and $\overline{f}$ induce the same action on $[A]^n$.

Let $f$ be in $W_2$. By
Definition~\ref{defn}, we get that $f=\beta_{n,n-1}^\ast g$, for some 
$g\in \cap_{i=1}^{k-1}(V_{A\setminus A_i}+\Ker\beta_{n,n-1}^\ast)$.
Lemma~\ref{lemma:aux} (applied to $k-1,n$ and $(A\setminus
A_1),\ldots,(A\setminus A_{k-1})$) yields
\begin{equation*}
\bigcap_{i=1}^{k-1}\left(V_{A\setminus
  A_i}+\Ker\beta_{n,n-1}^\ast\right)=\left(\bigcap_{i=1}^{k-1}V_{A\setminus
A_i}\right)+\Ker\beta_{n,n-1}^\ast.
\end{equation*}
Thence, up to replacing $g$ by $g+l$ (for some $l\in
\Ker\beta_{n,n-1}^\ast$), we may assume that $g\in
\cap_{i=1}^{k-1}V_{A\setminus   A_i}$. Let $\overline{g}$ be the
function in $\Fi_2^{[\Omega]^{n-1}}$ defined 
by
\[
\overline{g}(w)=\left\{
\begin{array}{ccl}
g(w)&&\textrm{if }w\subseteq A,\\
0&&\textrm{otherwise}.
\end{array}
\right.
\]
Set $\overline{f}=\beta_{n,n-1}^\ast \overline{g}$.  By construction,
$f$ and $\overline{f}$ coincide in $[A]^{n}$, that is,
$f$ and $\overline{f}$ induce the same action on $[A]^n$. 
Thus, it remains to prove that $\overline{f}\in W_1$, that is,
$\overline{f}$ vanishes on every $n$-subset $L$ of $(A\cap
A_i)\setminus A_i$, for $i=1,\ldots,k$. Let $L$ be an
$n$-subset of $(A\cup A_k)\setminus A_i$. We consider three cases
$L\subseteq A$, $|L\cap A_k|\geq 2$ and $|L\cap A_k|=1$.

If $L\subseteq A$, then
$\overline{f}(L)=f(L)=0$ (because $f$ and
$\overline{f}$ coincide on $[A]^n$).

If $|L\cap A_k|\geq 2$, then $(L\setminus \{x\})\nsubseteq A$, for every
$x$ in $L$. By definition of $\overline{g}$, we have 
$\overline{g}(L\setminus\{x\})=0$ and $\overline{f}(L)=\sum_{x\in
  L}\overline{g}(L\setminus \{x\})=0$.

If $|L\cap A_k|=1$ and $L\cap A_k=\{\overline{x}\}$, then (arguing as
in the previous paragraph)
$\overline{f}(L)=\sum_{x\in 
  L}\overline{g}(L\setminus\{x\})=g(L\setminus\{\overline{x}\})$. As
$L\subseteq (A\cup A_k)\setminus A_i$, we have that
$L\setminus\{\overline{x}\}\subseteq A\setminus A_i$. Since $g\in
V_{A\setminus A_i}$, we get that
$\overline{g}(L\setminus\{\overline{x}\})=g(L\setminus\{\overline{x}\})=0$. 
\cvd

J.Goodrick and A.Kolesnikov recently proved that if a complete stable theory
$T$ has $k$-uniqueness for every $2\leq k\leq n$, then $T$ has
$n+1$-existence~\cite{GoKo}. For completeness we report the proof of
their result.

\begin{theorem}\label{thm:GoKo}Let $T$ be a complete stable theory. If $T$ has
  $k$-uniqueness for 
  every $2\leq k\leq n$, then $T$ has $n+1$-existence.
\end{theorem}
\emph{Proof.}
Note that the existence and the uniqueness of nonforking extensions of
types in a stable theory yields that any stable theory has both
$2$-existence and $2$-uniqueness.

Since $T$ is a
complete stable theory, for every regular cardinal $k$, there exists a
saturated model of cardinality $k$. In the sequel we shall consider
the objects of $\mathcal{C}_T$ lying inside a very large saturated
``monster model'' $\mathfrak{C}$ of $T$.

Suppose $a$ is an $(n+1)$-amalgamation problem.  We have to prove that 
$a$ has  a solution $a'$. First, let $B_0$ and $B_1$ be sets of
$\mathfrak{C}$ such that
$\tp(B_0 / a(\emptyset)) 
= \tp(a([n])/ a(\emptyset))$, $\tp(B_1 /
a(\emptyset)) = \tp(a(\{n+1\})/ a(\emptyset))$, and $$B_0
\ind_{a(\emptyset)} B_1.$$  Let $\sigma_0$ and $\sigma_1$ be  two
automorphisms of $\mathfrak{C}$ fixing pointwise $a(\emptyset)$ and 
such 
that $B_0=\sigma_0(a([n]))$, $B_1=\sigma_1(a(\{n+1\}))$. 

Define $a'([n+1])$ to be the algebraic closure of
$B_0 \cup B_1$.  To determine the solution $a'$ of $a$, it remains to define 
the transition maps $a'_{s,[n+1]}:a'(s)\to a'([n+1])$, for all subsets
$s$ of 
$[n+1]$. The map $a'_{\emptyset,[n+1]}$ must be the
identity  on $a(\emptyset)$.  For $i$ in $[n]$, we let
$a'_{\{i\},[n+1]} : a(\{i\}) \to
a'([n+1])$ be the map $\sigma_0 \circ
a_{\{i\},[n]}$, and we let $a'_{\{n+1\},[n+1]}$ be the map $
\sigma_1$. Now, the following claim concludes the proof of the theorem.

\smallskip
\noindent\textsc{Claim: }
For every proper non-empty subset $s$ of $[n+1]$, there is a way to define the
transition maps 
$a'_{s,[n+1]}$, which is
consistent with $a$ and the definition of 
$a'_{\{i\},[n+1]}$ given above, and such that  $$a'_{s,[n+1]} (a(s)) =
\acl\left(\bigcup_{i \in s} a(\{i\})\right).$$ 

We argue by induction on the size $k$ of the set $s$. If $k=1$, then
there is nothing to prove.
Suppose we have defined $a'_{s,[n+1]}$ as in the claim, for all $s
\subseteq [n+1]$ such that $|s| < k$.  Let $s$ be a subset of 
$[n+1]$ such that $|s| = k$. The family of
sets $\{a(t)\mid t \subsetneq s \}$ forms a
$k$-amalgamation problem with the same transition maps as $a$. Call $a^1$
this amalgamation problem.  By the induction hypothesis, the
family of sets 
$\{a'_{t,[n+1]} (a(t)) \mid t \subsetneq s \}$ forms
another $k$-amalgamation problem with the transition maps
given by set inclusions. Call $a^2$ this amalgamation problem.  Notice
that $a^1$ and $a^2$ 
are isomorphic, and that both have independent solutions. Namely, $a^1$ can be
completed to $a(s)$ using the transition maps in $a$, and $a^2$ has a
natural solution $(a^2)'$ such that $$(a^2)'(s) = \acl\left(\bigcup_{i \in
  s} a(\{i\})\right),$$ where the transition maps are again given by set
inclusions.  So, by the $k$-uniqueness property, there is an
isomorphism of these solutions, which yields the desired transition
map $a'_{s,[n+1]}$ from $a(s)$ to $\acl(\bigcup_{i \in s} a(\{i\}))$. 

\cvd

Now we are ready to prove that $T_n$ has $k$-existence for every $k\leq n+1$.

\begin{prop}\label{exik}The theory $T_n$ has $k$-existence for every
  $k\leq n+1$.
\end{prop}
\emph{Proof.}
By definition, $T_n=\Th(M_n)$ is complete. Since $T_n$ is a stable
theory, the proof of this proposition follows at 
once from Proposition~\ref{unik} and Theorem~\ref{thm:GoKo}.
\cvd

Next, we show that $T_n$ does not have $n+1$-uniqueness.

\begin{prop}\label{unin+1}
The theory $T_n$ does not have $n+1$-uniqueness. 
 \end{prop}

\emph{Proof.}  Recall that by construction $n\geq 2$. Let
  $a:P(n+1)^-\to \mathcal{C}_{T_{n}}$ be the $(n+1)$-amalgamation problem
  defined on the objects by $a(s)=\acl(s)$ and where the morphisms are
  inclusions. In order to 
  prove this proposition we show the following equations: 
\begin{eqnarray}\label{eq1}
 |\Aut(\acl([n])/\cup_{i=1}^{n}\acl([n+1]-i))|&=&1, \\ 
 |\Aut(\acl([n])/\cup_{i=1}^{n}\acl([n]-i))|&=&2.\label{eq2}  
 \end{eqnarray}
In fact, by~\cite[Proposition~$3.5$]{Hr},
Equations~(\ref{eq1}),~(\ref{eq2}) yield that $a$ has more than one
solution up to isomorphism, i.e. $T_n$ does 
not have  $n+1$-uniqueness.  
 
We start by proving Equation~(\ref{eq1}). Since $[n],[n+1]-i$ have
size $n$, Proposition~\ref{acl1} yields 
 $\acl_{M_n}([n])=[n]\cup \{[n]\}\cup\{([n],0), ([n],1)\}$ and
  $\acl_{M_n}([n+1]-i)=([n+1]-i)\cup\{[n+1]-i\}\cup
\{([n+1]-i,0), ([n+1]-i,1)\}$. 

By the description given in the previous paragraph, every permutation
in $\Sym(\Omega)$ fixing 
pointwise the elements in $\cup_{i=1}^n\acl([n+1]-i)$ also fixes
pointwise every element in $\acl([n])$. Therefore, it suffices to
consider the elements in $\Imm\beta_{n,n-1}^\ast$. Let $f$ be in  
 $\Imm \beta^\ast_{n,n-1}$ and  suppose that $f$ fixes every element in
 $\cup_{i=1}^n\acl([n+1]-i)$, i.e. $f([n+1]-i)=0$, for $1\leq i\leq
 n$. Let $g\in\mathbb{F}_2^{[\Omega]^{n-1}}$ such that
 $f=\beta_{n,n-1}^\ast g$.  We have
 \begin{equation}\label{uniqueness}
0=\sum_{i=1}^{n} f([n+1]-i)=\sum_{i=1}^{n}\sum_{j\neq
  i}^{n+1}g([n+1]\setminus\{i,j\}). 
 \end{equation}
Now, for $j\neq n+1$, the summand $g([n+1]\setminus\{i,j\})$ appears
twice in Equation~(\ref{uniqueness}) and therefore over $\mathbb{F}_2$
their sum is zero. Hence
$$ 
0=\sum_{i=1}^{n} f([n+1]-i)=\sum _{i=1}^n g([n]-i)=(\beta_{n,n-1}^\ast
g)([n])=f([n]). 
$$
This yields that $f$ fixes $([n],0),([n],1)$. Hence
Equation~(\ref{eq1}) follows. 
 
We now prove Equation~(\ref{eq2}). Since $[n]-i$ has size $n-1$,
Proposition~\ref{acl1} implies $\acl_{M_n}([n]-i)=[n]-i$. Therefore,
$$\cup_{i=1}^n\acl_{M_n}([n]-i)=\cup_{i=1}^n([n]-i)=[n].$$ Also,
$\acl_{M_n}([n])=[n]\cup\{[n]\}\cup \{([n],0),([n],1)\}$. Corollary \ref{acl7} and Corollary \ref{acl5} yield that
every element of $\Aut(\acl([n])/\cup_{i=1}^{n}\acl([n]-i))$ fixes the
elements belonging to the sorts $\Omega$ and $[\Omega]^n$ of
$\acl_{M_n}([n])$. Hence, in order to prove
 Equation~(\ref{eq2}), it suffices to find an automorphism of
 $\acl_{M_n}([n])$ mapping $([n],0)$ into $([n],1)$. Let
 $g\in\Fi_2^{[\Omega]^{n-1}}$ with $g([n-1])=1$ and $g(w)=0$ for
 $w\neq [n-1]$. Set $f=\beta_{n,n-1}^\ast g$ and note that $f([n])=1$. As
 $\Aut(M_n)=\Imm\beta_{n,n-1}^\ast\rtimes \Sym(\Omega)$, the map $f$
 is an automorphism of $M_n$. By construction $f$ is an automorphism
 of $\acl_{M_n}([n])$ and $([n],0)^f=([n],0+f([n]))=([n],1)$.    
\cvd 

Finally, we show that $T_n$ does not have $n+2$-existence.
\begin{prop}\label{notn+2}
The theory $T_n$ does not have $n+2$-existence. 
\end{prop}
\emph{Proof.}
We construct an $n+2$-amalgamation problem $a$ over $\emptyset$ (that
is, $a(\emptyset)=\emptyset$) for $T_n$
with no solution.

Let $g$ be the element of $\mathbb{F}_2^{[\Omega]^{n-1}}$ defined by
\[
g(w)=\left\{
\begin{array}{ccl}
1&&\textrm{if }w=[n-1],\\
0&&\textrm{if }w\neq [n-1].
\end{array}
\right.
\]
Consider
$f=\beta_{n,n-1}^\ast g$ and note that, as
$\Aut(M_n)=\Imm\beta_{n,n-1}^\ast \rtimes \Sym(\Omega)$, the element
$f$ is an automorphism of $M_n$. 

Let $a$ be the functor
$a:P(n+2)^{-}\to \mathcal{C}_{T_n}$ defined on 
the objects by $a(s)=\acl(s)$ and with morphisms defined by

\begin{eqnarray}\label{eq:6}
a_{s,s'}=\left\{
\begin{array}{ccl}
f|_{a(s)}&&\textrm{if }s=[n] \textrm{ and } s'=[n+1],\\
\textrm{inclusion}&&\textrm{otherwise},
\end{array}
\right.
\end{eqnarray}
where $f|_{a(s)}$ denotes the restriction of the automorphism $f$  to
$a(s)$. It is not obvious from Equation~$(\ref{eq:6})$ 
that $a$ is a functor. Therefore, in the following paragraph, we  prove
that  $a$ is 
well-defined, that is, $a_{s_2,s_3}\circ a_{s_1,s_2}=a_{s_1,s_3}$ for every
$s_1,s_2,s_3$ in $P(n+2)^{-}$ with $s_1\subseteq s_2\subseteq s_3$. 

If $s_2\neq [n+1]$ and $s_3\neq [n+1]$, then (by Equation~$(\ref{eq:6})$)
the morphisms 
 $a_{s_1,s_2},a_{s_2,s_3}$ and 
$a_{s_1,s_3}$ are inclusions and so clearly $a_{s_2,s_3}\circ
a_{s_1,s_2}=a_{s_1,s_3}$. If $s_2=[n+1]$, then $s_2$ is a maximal
element of the partially ordered set $P(n+2)^{-}$. Thence
$s_3=s_2$ and, by Equation~$(\ref{eq:6})$,  
$a_{s_2,s_3}$ is the identity map. Thus $a_{s_2,s_3}\circ
a_{s_1,s_2}=a_{s_1,s_3}$. In particular, from now on we may assume that
$s_3=[n+1]$ and $s_2\neq [n+1]$. As $s_1\subseteq s_2$, if $s_2\neq
[n]$, then $s_1\neq [n]$ and so, by
Equation~$(\ref{eq:6})$, the morphisms 
 $a_{s_1,s_2},a_{s_2,s_3}$ and 
$a_{s_1,s_3}$ are inclusions and $a_{s_2,s_3}\circ
a_{s_1,s_2}=a_{s_1,s_3}$. If $s_2=s_1=[n]$, then $a_{s_1,s_2}$ is the
identity map and $a_{s_2,s_3}\circ
a_{s_1,s_2}=a_{s_1,s_3}$. The only case that remains to consider is
$s_3=[n+1]$, $s_2=[n]$ and $s_1\neq [n]$. Thence $a_{s_1,s_2}$ and
$a_{s_1,s_3}$ are inclusion maps and $a_{s_2,s_3}=f|_{a(s_2)}$. Since
$s_1\subseteq s_2=[n]$ and $s_1\neq [n]$, we have
$|s_1|<n-1$. Therefore,
$a(s_1)=\acl(s_1)$ and by Proposition~\ref{acl1} $\acl_{M_n}(s_1)=s_1$ consists only of elements belonging to the sort
$\Omega$ of $M_n$. As $f$ acts trivially on the elements belonging to
the sort $\Omega$, by Proposition \ref{acl2} we obtain $a_{s_2,s_3}\circ
a_{s_1,s_2}=(f|_{a(s_2)})|_{a(s_1)}=f|_{a(s_1)}=a_{s_1,s_3}$.
Finally, this proves that  $a:P(n+2)^{-}\to
\mathcal{C}_{T_n}$ is a functor.  

By Proposition~\ref{conclusione1},
$a(\emptyset)=\acl(\emptyset)=\emptyset$. Therefore, the functor $a$ is an
$n+2$-amalgamation problem over $\emptyset$ for $M_n$. 

We claim that
$a$ cannot be extended to 
$P(n+2)$. We argue by contradiction. Let $\overline{a}:P(n+2)\to
\mathcal{C}_{T_n}$ be a solution of $a$. In particular, $\overline{a}$
is an extension of $a$ to the whole of $P(n+2)$. 
Denote by $x_i$ the morphisms $\overline{a}_{[n+2]-i,[n+2]}$, for
$1\leq i\leq n+2$. So, by definition of morphism, $x_i$ is the
restriction to $\acl([n+2]-i)$  of an automorphism $f_i\sigma_i$ of
$M_n$, where $f_i\in\Imm\beta_{n,n-1}^\ast$ and  $\sigma_i\in \Sym(\Omega)$.

Since $\overline{a}$ is a functor and $\overline{a}$ extends $a$, we get
\begin{eqnarray}\label{eq4}
x_i\circ a_{[n+2]\setminus\{i,j\},[n+2]-i}&=&
\overline{a}_{[n+2]-i,[n+2]}\circ
\overline{a}_{[n+2]\setminus\{i,j\},[n+2]-i}\\\nonumber
&=&
\overline{a}_{[n+2]-j,[n+2]}\circ\overline{a}_{[n+2]\setminus\{i,j\},[n+2]-j}\\\nonumber  
&=&x_j\circ a_{[n+2]\setminus\{i,j\},[n+2]-j}.\nonumber
\end{eqnarray}
Let $i$ and $j$ be in $[n+2]$ with $i\neq j$. Fix an enumeration of $\acl_{M_n}([n+2]\setminus\{i,j\})$ and denote it as  $\overline{b_{ij}}=(b_{ij_1,},\dots)$. Then, as it is shown in Proposition \ref{acl2} $\overline{b_{ij}}\in \acl([n+2]\setminus\{i,j\})$ and, of course, also in $\acl([n+2]\setminus\{i\})$. By
Proposition~\ref{acl1} the ordered pair $([n+2]\setminus\{i,j\},0)$  belongs to the sort
$[\Omega]^n\times\Fi_2$ of $M_n$ and lies in $\acl_{M_n}([n+2]\setminus\{i,j\})$. Set $b_{ij_1}=([n+2]\setminus\{i,j\},0)$. We have  
\begin{eqnarray}\label{eq:1212}
x_i(\overline{b_{ij}})
&=&x_i(([n+2]\setminus\{i,j\},0),\dots)\\\nonumber
&=&((([n+2]\setminus\{i,j\})^{\sigma_i},0+f_i([n+2]\setminus\{i,j\})),\dots)\\\nonumber
&=&((([n+2]\setminus\{i,j\})^{\sigma_i},m_{ij}),\dots),\nonumber
\end{eqnarray} where

\begin{equation}\label{eq:1010}
m_{ij}=f_i([n+2]\setminus\{i,j\}).
\end{equation}
Consider
the matrix $M=(m_{ij})_{ij}$, with $m_{ii}=0$. 

Let $i$ and $j$ be in $[n+2]$ with $i\neq j$ and $\{i,j\}\neq
\{n+1,n+2\}$. By
Equation~$(\ref{eq:6})$ and by hypothesis on $\{i,j\}$, the morphism
$a_{[n+2]\setminus\{i,j\},[n+2]-i}$ is an inclusion map and so it fixes
$([n+2]\setminus\{i,j\},0)$. Therefore,
\begin{eqnarray*}
x_i\circ a_{[n+2]\setminus\{i,j\},[n+2]-i}(\overline{b_{ij}})
&=&x_i\circ a_{[n+2]\setminus\{i,j\},[n+2]-i}(([n+2]\setminus\{i,j\},0),\dots)\\
&=&x_i(([n+2]\setminus\{i,j\},0),\dots)\\
&=&((([n+2]\setminus\{i,j\})^{\sigma_i},m_{ij}),\dots),
\end{eqnarray*}
where in the last equality we used
Equations~$(\ref{eq:1212})$ and~$(\ref{eq:1010})$. Similarly,
replacing $i$ with $j$, we obtain 
\begin{eqnarray*}
x_i\circ a_{[n+2]\setminus\{i,j\},[n+2]-i}(\overline{b_{ij}})
&=&x_j\circ a_{[n+2]\setminus\{i,j\},[n+2]-j}(([n+2]\setminus\{i,j\},0),\dots)\\
&=&x_j(([n+2]\setminus\{i,j\},0),\dots)\\
&=&((([n+2]\setminus\{i,j\})^{\sigma_j},m_{ji}),\dots). 
\end{eqnarray*}
Now, by Equation~$(\ref{eq4})$, we have 

\begin{eqnarray*}x_i\circ a_{[n+2]\setminus\{i,j\},[n+2]-i}(\overline{b_{ij}})
&=&x_i\circ
a_{[n+2]\setminus\{i,j\},[n+2]-i}(([n+2]\setminus\{i,j\},0),\dots)\\
&=&x_j\circ
a_{[n+2]\setminus\{i,j\},[n+2]-j}(([n+2]\setminus\{i,j\},0),\dots).\end{eqnarray*}
In particular, 
\begin{equation}\label{eq:yyy}
m_{ij}=m_{ji}, \qquad \textrm{for every }i,j \textrm{ with }\{i,j\}\neq \{n+1,n+2\}.
\end{equation}

By Equation~$(\ref{eq:6})$ the morphism
$a_{[n+2]\setminus\{n+1,n+2\},[n+2]-(n+1)}$ is an inclusion map and so
it fixes
 $([n+2]\setminus\{n+1,n+2\},0)$. Therefore, 
\begin{eqnarray*}
& & x_{n+1}\circ a_{[n+2]\setminus\{n+1,n+2\},[n+2]-(n+1)}(\overline{b}_{n+1,n+2})\\
&=& x_{n+1}\circ a_{[n+2]\setminus\{n+1,n+2\},[n+2]-(n+1)}(([n+2]\setminus\{n+1,n+2\},0),\dots)\\
&=&x_{n+1}(([n+2]\setminus\{n+1,n+2\},0),\dots)\\
&=&((([n+2]\setminus\{n+1,n+2\})^{\sigma_{n+1}},m_{(n+1)(n+2)}),\dots).
\end{eqnarray*}
By Equation~$(\ref{eq:6})$ the morphism
$f|_{a([n])}=a_{[n],[n+1]}=a_{[n+2]\setminus\{n+1,n+2\},[n+2]-(n+2)}$
maps $([n+2]\setminus\{n+1,n+2\},0)$ to
$([n+2]\setminus\{n+1,n+2\},1)$. Therefore, 
\begin{eqnarray*}
& & x_{n+2}\circ a_{[n+2]\setminus\{n+1,n+2\},[n+2]-(n+2)}(\overline{b}_{n+1,n+2})\\
&=&x_{n+2}\circ a_{[n+2]\setminus\{n+1,n+2\},[n+2]-(n+2)}(([n+2]\setminus\{n+1,n+2\},0),\dots)\\
&=&x_{n+2}\circ f|_{a([n])}(([n+2]\setminus\{n+1,n+2\},0),\dots)\\
&=&x_{n+2}(([n+2]\setminus\{n+1,n+2\},1),\dots)\\
&=&((([n+2]\setminus\{n+1,n+2\})^{\sigma_{n+2}},m_{(n+2)(n+1)}+1),\dots).
\end{eqnarray*}
By Equation~$(\ref{eq4})$ (applied to $i=n+1$ and $j=n+2$), we have 
\begin{eqnarray*}
& & (([n+2]\setminus\{n+1,n+2\})^{\sigma_{n+1}},m_{(n+1)(n+2)})\\
& =& (([n+2]\setminus\{n+1,n+2\})^{\sigma_{n+2}},m_{(n+2)(n+1)}+1)
\end{eqnarray*}
and
\begin{equation}\label{eq:1414}
m_{(n+1)(n+2)}=m_{(n+2)(n+1)}+1.
\end{equation}

Now, we are ready to get a contradiction. We claim that  each row of
$M$ adds up to zero. We have
\begin{eqnarray*}
\sum_{j=1}^{n+2}m_{ij}&=&
\sum_{j\in ([n+2]-i)}m_{ij}=
\sum_{j\in ([n+2]-i)}f_i([n+2]\setminus\{i,j\})\\
&=&(\beta_{n+1,n}^\ast
f_i)([n+2]-i)=0, 
\end{eqnarray*}
where in the first equality we used that $m_{ii}=0$, in the second
equality we used Equation~$(\ref{eq:1010})$ and in the last equality
we used that $f_i\in
\Imm\beta_{n,n-1}^\ast=\Ker\beta_{n+1,n}^{*}$. In particular,
the sum of all the entries of 
$M$ is zero. Hence
$$0=\sum_{ij}m_{ij}=\sum_{i<j}(m_{ij}+m_{ji}).$$ By Equation~$(\ref{eq:yyy})$,
$m_{ij}=m_{ji}$ if $\{i,j\}\neq \{n+1,n+2\}$. So, in the previous sum there is
only one non-zero summand. Namely,
$m_{(n+1)(n+2)}+m_{(n+2)(n+1)}=0$. Now, Equation~$(\ref{eq:1414})$
yields $$m_{(n+1)(n+2)}+m_{(n+2)(n+1)}=m_{(n+1)(n+2)}+m_{(n+1)(n+2)}+1=1,$$ a
contradiction. This contradiction finally  
proves that the extension $\overline{a}$ does not exist.
\cvd

Now, Theorem~\ref{summary} follows at once from
Proposition~\ref{unik},~\ref{exik},~\ref{unin+1},~\ref{notn+2}. Finally,
we point out that Proposition~\ref{unin+1} also follows from
Theorem~\ref{thm:GoKo} and Proposition~\ref{notn+2}.

\section{Extension of Example~\ref{esempio1}}\label{coincidence}

In this section we remark that for every $n\geq 2$ the theories $T_n$ are stable and that the family of examples $\{M_n\}_{n\geq 2}$
generalizes the example due to E.Hrushovski given in~\cite{PKM}, see
Example~\ref{esempio1} in Section~\ref{sec0}.

\begin{defn}\label{generalization}
Let $\Omega$ be a countable set, and
    $C=[\Omega]^n\times \Z/2\Z$. Also let
    $E\subseteq \Omega\times [\Omega]^2$ be the membership relation, and let
    $P$ be the subset of $C^{n+1}$ such that
    $((w_1,\delta_1),\dots,(w_{n+1},\delta_{n+1}))\in P$ if and only if
    there are distinct $c_1,\dots,c_{n+1}\in \Omega$ such that
    $w_i=\{c_1,\dots, c_{n+1}\}\setminus c_i$ and 
    $\delta_1+\cdots+\delta_{n+1}=0$. Now let $\overline{M}_n$ be the model with the
    $3$-sorted universe $\Omega,[\Omega]^n,C$ and equipped with relations $E,P$
    and projection on the first coordinate $\pi:C\rightarrow
    [\Omega]^n$. Since $\overline{M}_n$ is a reduct of $(\Omega,
    \Z/2\Z)^{\textrm{eq}}$, we get that $\Th (\overline{M}_n)$ is stable. 
\end{defn}

\begin{prop}\label{coincides}
Let $\overline{M}_n$ be the structures described in Definition~\ref{generalization}. Then
$\Aut(\overline{M}_n)=\Imm\beta^\ast_{n,n-1} \rtimes \Sym(\Omega)$. In particular, $\overline{M}_n$
and $M_n$ are interdefinable.  
\end{prop}
\emph{Proof.} First we show that $\Sym(\Omega)$ is a subgroup  of
$\Aut(\overline{M}_n)$. Indeed, the group $\Sym(\Omega)$ acts with its natural action 
on the sorts $\Omega$ and $[\Omega]^n$ of $\overline{M}_n$. Also, 
if $g\in \Sym(\Omega)$ and $(\{a_1,\dots,a_n\},\delta)\in C$, then we set
$(\{a_1,\dots,a_n\},\delta)^g=(\{a_1^g,\dots,a_n^g\},\delta)$. This defines an
action of $\Sym(\Omega)$ on $\overline{M}_n$. It is
straightforward to see that the relations $E,P$ and the partition
given by the fibers of $\pi$  are preserved by
$\Sym(\Omega)$.
Hence, $\Sym(\Omega)\leq \Aut(\overline{M}_n)$. 

Let $\mu:\Aut(\overline{M}_n)\rightarrow \Sym(\Omega)$ be the  map given by
restriction on the sort $\Omega$ of $\overline{M}_n$. Since $\mu$ is a surjective
homomorphism, we 
have that $\Aut(\overline{M}_n)$ is a split extension of $\Ker\mu$ by
$\Sym(\Omega)$. Every element of $\Ker\mu$ preserves the fibres of
$\pi$  and fixes all the elements of 
$[\Omega]^n$. So $\Ker\mu$ is a closed $\Sym(\Omega)$-submodule of
$\Fi_2^{[\Omega]^n}$.

Let  $((w_1,\delta_1),\dots,(w_{n+1},\delta_{n+1}))$ be in $P$ and
$f$ be in
$\Ker\mu$. Since $\Ker \mu$ preserves $P$, we have $$f(w_1)+\delta_1
+\cdots +f(w_{n+1})+\delta_{n+1}=0.$$ 
From the definition of $P$ and $\beta_{n+1,n}^\ast$, we get
$$\Ker\mu=
\{f\in \Fi_2^{[\Omega]^n}\mid \sum_{x\in
  [w]^n}f(x)=0\textrm{ for every }w\in [\Omega]^{n+1}\}=\Ker\beta_{n+1,n}^\ast. $$
By Proposition~\ref{exact}, we have
that $\Ker\beta^\ast_{n+1,n}=\Imm\beta^\ast_{n,n-1}$. Therefore
$\Aut(\overline{M}_n)=\Aut(M_n)$ and $\overline{M}_n,M_n$ are interdefinable.
\cvd

\section*{Acknowledgements}

The authors thank J.Goodrick and A.Kolesnikov for
 the proof of Theorem~\ref{thm:GoKo} and for
giving their permission to include their proof in our paper. We are grateful to D. M. Evans for his stimulating suggestions and we thank
the anonymous referee for the very valuable comments and remarks on an earlier draft of the paper.

\end{document}